\documentclass[a4paper,12pt]{article}
\usepackage[utf8]{inputenc}
\usepackage[english]{babel}
\usepackage{amsmath,amssymb,amsfonts,amsthm}
\usepackage{tikz}
\usepackage{titling}
\usepackage{hyperref}
\usetikzlibrary{3d,calc}

\makeatletter
\def\@seccntformat#1{\csname the#1\endcsname.\ } 
\makeatother
\newtheorem{lemma}{Lemma}

\newtheorem{proposition}{Proposition}
\newcommand{\subtitle}[1]{%
  \posttitle{\par\end{center}\begin{center}\large#1\end{center}%
  \vskip 0.5em}%
}
\title{\bf \boldmath On the reconstruction of $n$-quasigroups of order $4$ and  upper bounds on their number\thanks{The work was supported by the RFBR grants 00-01-00822 and 99-01-00531.}}
\subtitle{(TRANSACTIONS OF THE CONFERENCE,
DEVOTED TO THE 90th ANNIVERSARY OF
ALEXEI A. LYAPUNOV,
Novosibirsk, October 8--12, 2001, pages 323--327)}
\author{D. S. Krotov, V. N. Potapov \\ Sobolev Institute of Mathematics, Novosibirsk}
\date{}

\def\FIG{

\begin{figure}[t]
$$
a)\!\!\!\!\!\!\!\!
\begin{tikzpicture}[scale=0.5,
    baseline={([yshift=-3mm] current bounding box.north)},
    x={(1cm,0cm)},
    y={(0.5cm,0.5cm)},
    z={(0cm,1.3cm)},
    dot/.style={circle, draw=black, inner sep=0pt, minimum size=3mm, line width=0.4pt}
]

\def\size{4}
\def\shift{1.7}

\foreach \z in {0,...,3} {
    \begin{scope}[shift={(0,0,\z * \shift)}]
        \foreach \y in {0,...,3} {
            \draw[gray!60, thin] (0,\y,0) -- (3,\y,0);
        }
        \foreach \x in {0,...,3} {
            \draw[gray!60, thin] (\x,0,0) -- (\x,3,0);
        }

        \foreach \x in {0,...,3} {
            \foreach \y in {0,...,3} {
                \pgfmathtruncatemacro{\iszero}{(\x==0 || \y==0 || \z==0) ? 1 : 0}
                \ifnum\iszero=1
                    \node[dot, fill=white] at (\x,\y,0) {};
                \else
                    \node[dot, fill=black] at (\x,\y,0) {};
                \fi
            }
        }
    \end{scope}
}
\draw[->, thick] (-1.3,-1.0,0) -- (1.0,-1.0,0);
\draw[->, thick] (-1.0,-1.3,0) -- (-1.0,1.0,0);
\draw[->, thick] (-1.0,-1.0,-0.1) -- (-1.0,-1.0,2.5);
\end{tikzpicture}
\qquad
b)
\begin{tikzpicture}[scale=0.5,
    baseline={([yshift=-3mm] current bounding box.north)},
    x={(1cm,0cm)},
    y={(0.5cm,0.5cm)},
    z={(0cm,1.3cm)},
    dot/.style={circle, draw=black, inner sep=0pt, minimum size=3mm, line width=0.4pt}
]

\def\size{4}
\def\shift{1.7}

\foreach \z in {0,...,3} {
    \begin{scope}[shift={(0,0,\z * \shift)}]
        \foreach \y in {0,...,3} {
            \draw[gray!60, thin] (0,\y,0) -- (3,\y,0);
        }
        \foreach \x in {0,...,3} {
            \draw[gray!60, thin] (\x,0,0) -- (\x,3,0);
        }

        \foreach \x in {0,...,3} {
            \foreach \y in {0,...,3} {
                \ifnum\z<2
                    \node[dot, fill=black] at (\x,\y,0) {};
                \else
                    \node[dot, fill=white] at (\x,\y,0) {};
                \fi
            }
        }
    \end{scope}
}
\end{tikzpicture}
\qquad
c)
\begin{tikzpicture}[scale=0.5,
    baseline={([yshift=-3mm] current bounding box.north)},
    x={(1cm,0cm)},
    y={(0.5cm,0.5cm)},
    z={(0cm,1.3cm)},
    dot/.style={circle, draw=black, inner sep=0pt, minimum size=3mm, line width=0.4pt}
]

\def\size{4}
\def\shift{1.7}

\foreach \z in {0,...,3} {
    \begin{scope}[shift={(0,0,\z * \shift)}]
        \foreach \y in {0,...,3} {
            \draw[gray!60, thin] (0,\y,0) -- (3,\y,0);
        }
        \foreach \x in {0,...,3} {
            \draw[gray!60, thin] (\x,0,0) -- (\x,3,0);
        }

        \foreach \x in {0,...,3} {
            \foreach \y in {0,...,3} {
                \pgfmathtruncatemacro{\iszero}{(\x==0 || \y==0 || \z==0) ? 1 : 0}
                \ifnum\iszero=1
                    \node[dot, fill=black] at (\x,\y,0) {};
                \else
                    \node[dot, fill=white] at (\x,\y,0) {};
                \fi
            }
        }
    \end{scope}
}
\end{tikzpicture}
\qquad
d)
\begin{tikzpicture}[scale=0.5,
    baseline={([yshift=-3mm] current bounding box.north)},
    x={(1cm,0cm)},
    y={(0.5cm,0.5cm)},
    z={(0cm,1.3cm)},
    dot/.style={circle, draw=black, inner sep=0pt, minimum size=3mm, line width=0.4pt}
]

\def\size{4}
\def\shift{1.7}

\foreach \z in {0,...,3} {
    \begin{scope}[shift={(0,0,\z * \shift)}]
        \foreach \y in {0,...,3} {
            \draw[gray!60, thin] (0,\y,0) -- (3,\y,0);
        }
        \foreach \x in {0,...,3} {
            \draw[gray!60, thin] (\x,0,0) -- (\x,3,0);
        }

        \foreach \x in {0,...,3} {
            \foreach \y in {0,...,3} {
                \pgfmathtruncatemacro{\numzeros}{(\x==0 ? 1 : 0) + (\y==0 ? 1 : 0) + (\z==0 ? 1 : 0)}
                \ifnum\numzeros>1
                    \node[dot, fill=black] at (\x,\y,0) {};
                \else
                    \node[dot, fill=white] at (\x,\y,0) {};
                \fi
            }
        }
    \end{scope}
}
\end{tikzpicture}
$$
\caption{Subsets of $\{0,1,2,3\}^3$: a) $G'_3$, b) $G''_3$, c) $G'''_3$, d) $G^1_3$.}\label{fig}
\end{figure}}

\begin{document}

\maketitle
\vspace{-3ex}
\begin{abstract}
The number of ways for partial $n$-quasigroups of order $4$ to be extended to $n$-quasigroups is investigated. Upper bounds on the number of $n$-quasigroups of order $4$ are obtained. The asymptotic $3^{n+1}2^{2^n+1}$ of that number is established.
\end{abstract}

Let $n$ be an arbitrary natural number, $A$ a finite%
\footnote{For the case of an infinite $A$, the definition of $n$-quasigroup should be given in another way.} set of cardinality $m>0$, and $G \subseteq A^n$. Let $q$ be a mapping from $G$ to $A$. If $q(c) \neq q(c')$ holds for every $c,c' \in G$ that differ in exactly one coordinate then the triple $(A,G,q)$ is called a \emph{partial $n$-quasigroup} (of order~$m$). If in addition $G = A^n$ then the pair $(A,q)$ is called an \emph{$n$-quasigroup} (of order~$m$).

It is natural to represent an $n$-quasigroup $(A,q)$
by the array of values of $q$. In such an $n$-dimensional array of size
$|A|\times\cdots\times|A|$, every line in each coordinate
direction contains all the elements of $A$
(see, for example, Fig.~\ref{fig2}).

Two $n$-quasigroups $(A,q)$ and $(A,\hat{q})$ are called \emph{equivalent} if there exist permutations $\rho_0, \rho_1, \ldots, \rho_n : A \to A$ and $\pi : \{1,\ldots,n\} \to \{1,\ldots,n\}$ such that
\[
\hat{q}(x_1,\ldots,x_n) = \rho_0 q(\rho_1 x_{\pi(1)}, \ldots, \rho_n x_{\pi(n)}).
\]

The $n$-quasigroup $(A,q)$ is called an \emph{extension} of the partial $n$-quasigroup $(A,G,p)$ if $p = q|_G$. Note that a partial $n$-quasigroup may have more than one extension or no extension at all.

There are exactly two $n$-quasigroups of order $2$. All $n$-quasigroups of order 3 are equivalent, and their number is $3 \cdot 2^n$. The case of order $4$ is the first nontrivial one
from the viewpoint of diversity.
Below, the order of all (partial) quasigroups is~$4$ by default.

Let $A = \{0,1,2,3\}$ and let $Q(n)$ be the number of different $n$-quasigroups $(A,q)$. It was shown in~\cite{Krotov2000} that
\begin{equation}\label{eq1}
Q(n) \geq L(n) \stackrel{\text{df}}{=} 3^{n+1}2^{2^n+1} - 2^{n+3}3^n.
\end{equation}

The upper bounds \eqref{eq2}--\eqref{eq6} on $Q(n)$ are based on estimates
for the numbers of extensions of certain partial $n$-quasigroups.

\section{The trivial bound}
\FIG
Let $G_n' = \{1,2,3\}^n$ (see, e.g., Fig.~\ref{fig}a), and let $(A,G_n',p)$ be a partial $n$-quasigroup.

\begin{lemma}\label{l1}
A partial $n$-quasigroup $(A,G_n',p)$ has at most one extension.
\end{lemma}

Since there are $4^{3^n}$ maps from $G_n'$ to $A$, we have
\begin{equation}\label{eq2}
Q(n) \leq 4^{3^n} = 2^{2 ^{\log_23 \cdot n+1}}.
\end{equation}

\section{Asymptotics of \texorpdfstring{$\log\log Q(n)$}{log log Q(n)}}
\label{s2}
Let $G_n'' = A^{n-1} \times \{0,1\}$  (see, e.g., Fig.~\ref{fig}b), and let $(A,G_n'',p)$ be a partial $n$-quasigroup.

\begin{lemma}\label{l2}
The partial $n$-quasigroup $(A,G_n'',p)$ has at most $2^{2^{n-1}}$ extensions.
\end{lemma}

The mapping $p$ is uniquely defined by two subfunctions
\[
p^{(0)}(x_1,\ldots,x_{n-1}) \stackrel{\text{df}}{=} p(x_1,\ldots,x_{n-1},0)
\]
and
\[
p^{(1)}(x_1,\ldots,x_{n-1}) \stackrel{\text{df}}{=} p(x_1,\ldots,x_{n-1},1).
\]
The pairs $(A,p^{(0)})$ and $(A,p^{(1)})$ are $(n-1)$-quasigroups. Therefore Lemma~\ref{l2} implies
\[
Q(n) \leq 2^{2^{n-1}} (Q(n-1))^2
\]
and, by induction,
\begin{equation}\label{eq3}
Q(n) \leq 2^{(n+\text{const})2^{n-1}} = 2^{2^{n+\log_2 n - 1 + o(1)}}.
\end{equation}
Bounds~\eqref{eq1} and~\eqref{eq3} give the asymptotic of $\log\log Q(n)$:
\[
\log_2 \log_2 Q(n) = n + O(\log_2 n).
\]

\section{Asymptotics of \texorpdfstring{$\log Q(n)$}{log Q(n)}}

Let $G'''_n = A^n \setminus \{1,2,3\}^n$ (see, e.g., Fig.~\ref{fig}c), and let $(A,G'''_n,p)$ be a partial $n$-quasigroup.

\begin{lemma}\label{l3}
If $n \geq 2$ then the partial $n$-quasigroup $(A,G'''_n,p)$ has at most four extensions.
\end{lemma}

For an $n$-ary mapping $q$, we call a $k$-ary \emph{minorant} the mapping obtained by substituting zeros for some $n-k$ arguments of~$q$. Lemma~\ref{l3} implies that a $k$-ary minorant of~$q$ can be reconstructed from the $(k-1)$-ary minorants in at most four ways, provided $(A,q)$ is an $n$-quasigroup and $k>1$.

Let $G_n^1 = \bigcup_{i=1}^n A_i$, where
\(
A_i = \{x = (x_1,\ldots,x_n) \in A^n : x_j = 0 \text{ if } j \neq i\}
\)  (see, e.g., Fig.~\ref{fig}d).
Let $(A,q)$ be an $n$-quasigroup. The partial $n$-quasigroup $(A,G_n^1,p)$, where $p = q|_{G_n^1}$, specifies all the $1$-ary minorants of~$q$. There are $2^n$ minorants in total, of which $n+1$ are $0$- and $1$-ary. Therefore, assuming $p$ is specified, the number of ways for $q$ to be reconstructed is at most $4^{2^n - n - 1}$. The number of ways to specify $p$ is $4 \cdot 6^n$, therefore
\begin{equation}\label{eq4}
Q(n) \leq 4 \cdot 6^n \cdot 4^{2^n - n - 1} = 2^{2^{n+1} + O(n)}.
\end{equation}

Let the mapping $\oplus : A^2 \to A$ be defined by the value table
\[
\begin{array}{c|cccc}
\oplus & 0 & 1 & 2 & 3 \\
\hline
0 & 0 & 1 & 2 & 3 \\
1 & 1 & 0 & 3 & 2 \\
2 & 2 & 3 & 1 & 0 \\
3 & 3 & 2 & 0 & 1
\end{array}
\]
Let the mapping $b : A \to \{0,1\}$ be defined by the identities
\[
b(0) = b(1) = 0, \qquad b(2) = b(3) = 1.
\]
Call an $n$-quasigroup $(A,q)$ \emph{bicubical} if it is equivalent to the $n$-quasigroup $(A,q'')$ with
\[
q''(x_1,\ldots,x_n) = x_1 \oplus \cdots \oplus x_n.
\]
Call an  $n$-quasigroup $(A,q)$ \emph{cubical}%
\footnote{In the further study of $n$-ary quasigroups of order $4$,
the term ``bicubical'' is replaced by ``linear'' and ``cubical'' by ``semilinear'' (reprint remark).}
if it is equivalent to some $n$-quasigroup $(A,q')$ such that
\[
b(q'(x_1,\ldots,x_n)) = b(x_1 \oplus \cdots \oplus x_n).
\]
Fig.~\ref{fig2} gives examples of bicubical, cubical (but not bicubical), and noncubical $3$-quasigroups.

\begin{proposition}\label{p1}
The number of cubical $n$-quasigroups is $L(n)$ (see~\eqref{eq1}).
\end{proposition}

\begin{lemma}\label{l4}
Let $n \geq 2$, and $(A,G'''_n,p)$ be a partial $n$-quasigroup. Then
\begin{enumerate}
\item if $(A,G'''_n,p)$ has more than one extension then all the extensions are cubical;
\item if $(A,G'''_n,p)$ has more than two extensions then it has a bicubical extension.
\end{enumerate}
\end{lemma}

Lemma~\ref{l4} yields the following improvement of the preceding bound:
\begin{equation}\label{eq5}
Q(n) \leq 4 \cdot 6^n \cdot 3^{\binom{n}{\lfloor n/2\rfloor}} \cdot 2^{2^n} = 2^{2^n (1 + O(1/\sqrt{n}))}.
\end{equation}

It follows from \eqref{eq1} and \eqref{eq5} that
\[
\log_2 Q(n) = 2^n (1 + O(1/\sqrt{n})).
\]

\section{Asymptotics of \texorpdfstring{$Q(n)$}{Q(n)}}

The $n$-quasigroup $(A, q)$ is called \emph{decomposable} (permutably reducible) if there exist $k \geq 2$, a $k$-quasigroup $(A, \dot{q})$, an $(n-k+1)$-quasigroup $(A, \ddot{q})$, and a permutation $\pi : \{1, \ldots, n\} \to \{1, \ldots, n\}$ such that
\[
q(x_1, \ldots, x_n) = \dot{q}(x_{\pi(1)}, \ldots, x_{\pi(k-1)}, \ddot{q}(x_{\pi(k)}, \ldots, x_{\pi(n)})).
\]
Fig.~\ref{fig2}c gives an example of a decomposable $3$-quasigroup.

\begin{figure}[t]
\scalebox{0.8}{
$
\mbox{}\qquad\qquad
\text{\Large a)}\!\!\!\!\!\!\!\!\!\!\!\!\!\!\!\!\!\!\!\!\!\!\!\!\!\!\!
\begin{tikzpicture}[
    baseline={([yshift=-3mm] current bounding box.north)},
    x={(1cm,0cm)},
    y={(0.5cm,0.5cm)},
    z={(0cm,1.3cm)},
    dot/.style={circle, draw=black, inner sep=0pt, minimum size=5mm, line width=0.4pt, font=\small},
    dotbold/.style={circle, draw=black, inner sep=0pt, minimum size=5mm, line width=1.5pt, font=\small}
]

\def\size{4}
\def\shift{1.7}

\foreach \z in {0,...,3} {
    \begin{scope}[shift={(0,0,\z * \shift)}]
        \foreach \y in {0,...,3} {
            \draw[gray!60, thin] (0,\y,0) -- (3,\y,0);
        }
        \foreach \x in {0,...,3} {
            \draw[gray!60, thin] (\x,0,0) -- (\x,3,0);
        }

        \ifnum\z=0
            \node[dotbold, fill=white] at (0,0,0) {\bf 0};
            \node[dotbold, fill=white] at (1,0,0) {\bf 1};
            \node[dot, fill=white] at (2,0,0) {\bf 2};
            \node[dot, fill=white] at (3,0,0) {\bf 3};
            \node[dotbold, fill=white] at (0,1,0) {\bf 1};
            \node[dotbold, fill=white] at (1,1,0) {\bf 0};
            \node[dot, fill=white] at (2,1,0) {\bf 3};
            \node[dot, fill=white] at (3,1,0) {\bf 2};
            \node[dot, fill=white] at (0,2,0) {\bf 2};
            \node[dot, fill=white] at (1,2,0) {\bf 3};
            \node[dot, fill=white] at (2,2,0) {\bf 0};
            \node[dot, fill=white] at (3,2,0) {\bf 1};
            \node[dot, fill=white] at (0,3,0) {\bf 3};
            \node[dot, fill=white] at (1,3,0) {\bf 2};
            \node[dot, fill=white] at (2,3,0) {\bf 1};
            \node[dot, fill=white] at (3,3,0) {\bf 0};
        \fi

        \ifnum\z=1
            \node[dotbold, fill=white] at (0,0,0) {\bf 1};
            \node[dotbold, fill=white] at (1,0,0) {\bf 0};
            \node[dot, fill=white] at (2,0,0) {\bf 3};
            \node[dot, fill=white] at (3,0,0) {\bf 2};
            \node[dotbold, fill=white] at (0,1,0) {\bf 0};
            \node[dotbold, fill=white] at (1,1,0) {\bf 1};
            \node[dot, fill=white] at (2,1,0) {\bf 2};
            \node[dot, fill=white] at (3,1,0) {\bf 3};
            \node[dot, fill=white] at (0,2,0) {\bf 3};
            \node[dot, fill=white] at (1,2,0) {\bf 2};
            \node[dot, fill=white] at (2,2,0) {\bf 1};
            \node[dot, fill=white] at (3,2,0) {\bf 0};
            \node[dot, fill=white] at (0,3,0) {\bf 2};
            \node[dot, fill=white] at (1,3,0) {\bf 3};
            \node[dot, fill=white] at (2,3,0) {\bf 0};
            \node[dot, fill=white] at (3,3,0) {\bf 1};
        \fi

        \ifnum\z=2
            \node[dot, fill=white] at (0,0,0) {\bf 2};
            \node[dot, fill=white] at (1,0,0) {\bf 3};
            \node[dot, fill=white] at (2,0,0) {\bf 0};
            \node[dot, fill=white] at (3,0,0) {\bf 1};
            \node[dot, fill=white] at (0,1,0) {\bf 3};
            \node[dot, fill=white] at (1,1,0) {\bf 2};
            \node[dot, fill=white] at (2,1,0) {\bf 1};
            \node[dot, fill=white] at (3,1,0) {\bf 0};
            \node[dot, fill=white] at (0,2,0) {\bf 0};
            \node[dot, fill=white] at (1,2,0) {\bf 1};
            \node[dot, fill=white] at (2,2,0) {\bf 2};
            \node[dot, fill=white] at (3,2,0) {\bf 3};
            \node[dot, fill=white] at (0,3,0) {\bf 1};
            \node[dot, fill=white] at (1,3,0) {\bf 0};
            \node[dot, fill=white] at (2,3,0) {\bf 3};
            \node[dot, fill=white] at (3,3,0) {\bf 2};
        \fi

        \ifnum\z=3
            \node[dot, fill=white] at (0,0,0) {\bf 3};
            \node[dot, fill=white] at (1,0,0) {\bf 2};
            \node[dot, fill=white] at (2,0,0) {\bf 1};
            \node[dot, fill=white] at (3,0,0) {\bf 0};
            \node[dot, fill=white] at (0,1,0) {\bf 2};
            \node[dot, fill=white] at (1,1,0) {\bf 3};
            \node[dot, fill=white] at (2,1,0) {\bf 0};
            \node[dot, fill=white] at (3,1,0) {\bf 1};
            \node[dot, fill=white] at (0,2,0) {\bf 1};
            \node[dot, fill=white] at (1,2,0) {\bf 0};
            \node[dot, fill=white] at (2,2,0) {\bf 3};
            \node[dot, fill=white] at (3,2,0) {\bf 2};
            \node[dot, fill=white] at (0,3,0) {\bf 0};
            \node[dot, fill=white] at (1,3,0) {\bf 1};
            \node[dot, fill=white] at (2,3,0) {\bf 2};
            \node[dot, fill=white] at (3,3,0) {\bf 3};
        \fi
    \end{scope}
}
\draw[->, thick] (-1.3,-1.0,0) -- (1.0,-1.0,0);
\draw[->, thick] (-1.0,-1.3,0) -- (-1.0,1.0,0);
\draw[->, thick] (-1.0,-1.0,-0.1) -- (-1.0,-1.0,2.5);
\end{tikzpicture}
\quad
\qquad
\text{\Large b)}\!\!\!\!
\begin{tikzpicture}[
    baseline={([yshift=-3mm] current bounding box.north)},
    x={(1cm,0cm)},
    y={(0.5cm,0.5cm)},
    z={(0cm,1.3cm)},
    dot/.style={circle, draw=black, inner sep=0pt, minimum size=5mm, line width=0.4pt, font=\small, fill=white},
    dotbold/.style={circle, draw=black, inner sep=0pt, minimum size=5mm, line width=1.5pt, font=\small, fill=white}
]

\def\size{4}
\def\shift{1.7}
\def\radius{2.5mm}  


\draw[<->, thick, >=stealth, shorten <= \radius, shorten >= \radius] (2,2,0) -- (2,2,1.7);
\draw[<->, thick, >=stealth, shorten <= \radius, shorten >= \radius] (3,2,0) -- (3,2,1.7);
\draw[<->, thick, >=stealth, shorten <= \radius, shorten >= \radius] (2,3,0) -- (2,3,1.7);
\draw[<->, thick, >=stealth, shorten <= \radius, shorten >= \radius] (3,3,0) -- (3,3,1.7);
\foreach \z in {0,...,3} {
    \begin{scope}[shift={(0,0,\z * \shift)}]
        \foreach \y in {0,...,3} {
            \draw[gray!60, thin] (0,\y,0) -- (3,\y,0);
        }
        \foreach \x in {0,...,3} {
            \draw[gray!60, thin] (\x,0,0) -- (\x,3,0);
        }

        \ifnum\z=0
            \draw[<->, thick, >=stealth, shorten <= \radius, shorten >= \radius] (2,2,0) -- (3,2,0);
            \draw[<->, thick, >=stealth, shorten <= \radius, shorten >= \radius] (2,3,0) -- (3,3,0);
            \draw[<->, thick, >=stealth, shorten <= \radius, shorten >= \radius] (2,2,0) -- (2,3,0);
            \draw[<->, thick, >=stealth, shorten <= \radius, shorten >= \radius] (3,2,0) -- (3,3,0);
        \fi

        \ifnum\z=1
            \draw[<->, thick, >=stealth, shorten <= \radius, shorten >= \radius] (2,2,0) -- (3,2,0);
            \draw[<->, thick, >=stealth, shorten <= \radius, shorten >= \radius] (2,3,0) -- (3,3,0);
            \draw[<->, thick, >=stealth, shorten <= \radius, shorten >= \radius] (2,2,0) -- (2,3,0);
            \draw[<->, thick, >=stealth, shorten <= \radius, shorten >= \radius] (3,2,0) -- (3,3,0);
        \fi

        \ifnum\z=0
            \node[dot] at (0,0,0) {\bf 0};
            \node[dot] at (1,0,0) {\bf 1};
            \node[dot] at (2,0,0) {\bf 2};
            \node[dot] at (3,0,0) {\bf 3};
            \node[dot] at (0,1,0) {\bf 1};
            \node[dot] at (1,1,0) {\bf 0};
            \node[dot] at (2,1,0) {\bf 3};
            \node[dot] at (3,1,0) {\bf 2};
            \node[dot] at (0,2,0) {\bf 2};
            \node[dot] at (1,2,0) {\bf 3};
            \node[dotbold] at (2,2,0) {\bf 1};
            \node[dotbold] at (3,2,0) {\bf 0};
            \node[dot] at (0,3,0) {\bf 3};
            \node[dot] at (1,3,0) {\bf 2};
            \node[dotbold] at (2,3,0) {\bf 0};
            \node[dotbold] at (3,3,0) {\bf 1};
        \fi

        \ifnum\z=1
            \node[dot] at (0,0,0) {\bf 1};
            \node[dot] at (1,0,0) {\bf 0};
            \node[dot] at (2,0,0) {\bf 3};
            \node[dot] at (3,0,0) {\bf 2};
            \node[dot] at (0,1,0) {\bf 0};
            \node[dot] at (1,1,0) {\bf 1};
            \node[dot] at (2,1,0) {\bf 2};
            \node[dot] at (3,1,0) {\bf 3};
            \node[dot] at (0,2,0) {\bf 3};
            \node[dot] at (1,2,0) {\bf 2};
            \node[dotbold] at (2,2,0) {\bf 0};
            \node[dotbold] at (3,2,0) {\bf 1};
            \node[dot] at (0,3,0) {\bf 2};
            \node[dot] at (1,3,0) {\bf 3};
            \node[dotbold] at (2,3,0) {\bf 1};
            \node[dotbold] at (3,3,0) {\bf 0};
        \fi

        \ifnum\z=2
            \node[dot] at (0,0,0) {\bf 2};
            \node[dot] at (1,0,0) {\bf 3};
            \node[dot] at (2,0,0) {\bf 0};
            \node[dot] at (3,0,0) {\bf 1};
            \node[dot] at (0,1,0) {\bf 3};
            \node[dot] at (1,1,0) {\bf 2};
            \node[dot] at (2,1,0) {\bf 1};
            \node[dot] at (3,1,0) {\bf 0};
            \node[dot] at (0,2,0) {\bf 0};
            \node[dot] at (1,2,0) {\bf 1};
            \node[dot] at (2,2,0) {\bf 2};
            \node[dot] at (3,2,0) {\bf 3};
            \node[dot] at (0,3,0) {\bf 1};
            \node[dot] at (1,3,0) {\bf 0};
            \node[dot] at (2,3,0) {\bf 3};
            \node[dot] at (3,3,0) {\bf 2};
        \fi

        \ifnum\z=3
            \node[dot] at (0,0,0) {\bf 3};
            \node[dot] at (1,0,0) {\bf 2};
            \node[dot] at (2,0,0) {\bf 1};
            \node[dot] at (3,0,0) {\bf 0};
            \node[dot] at (0,1,0) {\bf 2};
            \node[dot] at (1,1,0) {\bf 3};
            \node[dot] at (2,1,0) {\bf 0};
            \node[dot] at (3,1,0) {\bf 1};
            \node[dot] at (0,2,0) {\bf 1};
            \node[dot] at (1,2,0) {\bf 0};
            \node[dot] at (2,2,0) {\bf 3};
            \node[dot] at (3,2,0) {\bf 2};
            \node[dot] at (0,3,0) {\bf 0};
            \node[dot] at (1,3,0) {\bf 1};
            \node[dot] at (2,3,0) {\bf 2};
            \node[dot] at (3,3,0) {\bf 3};
        \fi
    \end{scope}
}
\end{tikzpicture}
\quad
\qquad
\text{\Large c)}\!\!\!\!
\begin{tikzpicture}[
    baseline={([yshift=-3mm] current bounding box.north)},
    x={(1cm,0cm)},
    y={(0.5cm,0.5cm)},
    z={(0cm,1.3cm)},
    dot/.style={circle, draw=black, inner sep=0pt, minimum size=5mm, line width=0.4pt, font=\small, fill=white},
    dotbold/.style={circle, draw=black, inner sep=0pt, minimum size=5mm, line width=1.5pt, font=\small, fill=white},
    dotgray/.style={circle, draw=black, inner sep=0pt, minimum size=5mm, line width=0.4pt, font=\small, fill=lightgray},
    dotboldgray/.style={circle, draw=black, inner sep=0pt, minimum size=5mm, line width=1.5pt, font=\small, fill=lightgray}
]

\def\size{4}
\def\shift{1.7}
\def\radius{2.5mm}

\foreach \z in {0,...,3} {
    \begin{scope}[shift={(0,0,\z * \shift)}]
        \foreach \y in {0,...,3} {
            \draw[gray!60, thin] (0,\y,0) -- (3,\y,0);
        }
        \foreach \x in {0,...,3} {
            \draw[gray!60, thin] (\x,0,0) -- (\x,3,0);
        }

        \ifnum\z=0
            \node[dotboldgray] at (0,0,0) {\bf 0};
            \node[dotbold] at (1,0,0) {\bf 1};
            \node[dotbold] at (2,0,0) {\bf 2};
            \node[dotbold] at (3,0,0) {\bf 3};
            \node[dotboldgray] at (0,1,0) {\bf 1};
            \node[dotbold] at (1,1,0) {\bf 3};
            \node[dotbold] at (2,1,0) {\bf 0};
            \node[dotbold] at (3,1,0) {\bf 2};
            \node[dotboldgray] at (0,2,0) {\bf 2};
            \node[dotbold] at (1,2,0) {\bf 0};
            \node[dotbold] at (2,2,0) {\bf 3};
            \node[dotbold] at (3,2,0) {\bf 1};
            \node[dotboldgray] at (0,3,0) {\bf 3};
            \node[dotbold] at (1,3,0) {\bf 2};
            \node[dotbold] at (2,3,0) {\bf 1};
            \node[dotbold] at (3,3,0) {\bf 0};
        \fi

        \ifnum\z=1
            \node[dotgray] at (0,0,0) {\bf 1};
            \node[dot] at (1,0,0) {\bf 0};
            \node[dot] at (2,0,0) {\bf 3};
            \node[dot] at (3,0,0) {\bf 2};
            \node[dotgray] at (0,1,0) {\bf 0};
            \node[dot] at (1,1,0) {\bf 2};
            \node[dot] at (2,1,0) {\bf 1};
            \node[dot] at (3,1,0) {\bf 3};
            \node[dotgray] at (0,2,0) {\bf 3};
            \node[dot] at (1,2,0) {\bf 1};
            \node[dot] at (2,2,0) {\bf 2};
            \node[dot] at (3,2,0) {\bf 0};
            \node[dotgray] at (0,3,0) {\bf 2};
            \node[dot] at (1,3,0) {\bf 3};
            \node[dot] at (2,3,0) {\bf 0};
            \node[dot] at (3,3,0) {\bf 1};
        \fi

        \ifnum\z=2
            \node[dotgray] at (0,0,0) {\bf 2};
            \node[dot] at (1,0,0) {\bf 3};
            \node[dot] at (2,0,0) {\bf 1};
            \node[dot] at (3,0,0) {\bf 0};
            \node[dotgray] at (0,1,0) {\bf 3};
            \node[dot] at (1,1,0) {\bf 0};
            \node[dot] at (2,1,0) {\bf 2};
            \node[dot] at (3,1,0) {\bf 1};
            \node[dotgray] at (0,2,0) {\bf 1};
            \node[dot] at (1,2,0) {\bf 2};
            \node[dot] at (2,2,0) {\bf 0};
            \node[dot] at (3,2,0) {\bf 3};
            \node[dotgray] at (0,3,0) {\bf 0};
            \node[dot] at (1,3,0) {\bf 1};
            \node[dot] at (2,3,0) {\bf 3};
            \node[dot] at (3,3,0) {\bf 2};
        \fi

        \ifnum\z=3
            \node[dotgray] at (0,0,0) {\bf 3};
            \node[dot] at (1,0,0) {\bf 2};
            \node[dot] at (2,0,0) {\bf 0};
            \node[dot] at (3,0,0) {\bf 1};
            \node[dotgray] at (0,1,0) {\bf 2};
            \node[dot] at (1,1,0) {\bf 1};
            \node[dot] at (2,1,0) {\bf 3};
            \node[dot] at (3,1,0) {\bf 0};
            \node[dotgray] at (0,2,0) {\bf 0};
            \node[dot] at (1,2,0) {\bf 3};
            \node[dot] at (2,2,0) {\bf 1};
            \node[dot] at (3,2,0) {\bf 2};
            \node[dotgray] at (0,3,0) {\bf 1};
            \node[dot] at (1,3,0) {\bf 0};
            \node[dot] at (2,3,0) {\bf 2};
            \node[dot] at (3,3,0) {\bf 3};
        \fi
    \end{scope}
}
\end{tikzpicture}
$
}
\caption{Examples of $3$-quasigroups:\\
    a) bicubical $3$-quasigroup $(A, q'')$;
    the value array can be partitioned in three ways into
$2\times2\times2$ subarrays, each containing exactly two distinct elements;\\
    b) cubical $3$-quasigroup obtained from $(A, q'')$ by ``switching'' one $2 \times 2 \times 2$-subarray;\\
    c) decomposable $3$-quasigroup $(A, q)$ which is not cubical; $q(x, y, z) = q_{\text{out}}(q_{\text{in}}(x, y), z)$; in the figure the value arrays of $q_{\text{out}}$ and $q_{\text{in}}$ are marked by gray and bold respectively.}\label{fig2}
\end{figure}

Let $D_n$ be the set of all decomposable $n$-quasigroups $(A, q)$, and let $\overline C_n$ be the set of $n$-quasigroups $(A, q)$ which are not cubical.

Inequality~\eqref{eq5} readily implies the following proposition.

\begin{proposition}\label{p2}
If $n \geq 6$ then $|D_n| \leq 2^{2^n}$.
\end{proposition}

The following statement is established by an exhaustive computer enumeration of $5$-quasigroups of order~$4$.

\begin{proposition}\label{p3}
$|\overline C_5| = 4 \cdot 6^5 \cdot 201538000 - L(5) = 4 \cdot 6^5 \cdot 211410 \leq 2 \cdot 2^{2^5}$.
\end{proposition}

Let $G''_n$, $p^{(0)}$, and $p^{(1)}$ be defined as in Section~\ref{s2}.

\begin{lemma}\label{l5}
Let $(A, q) \in \overline C_n \setminus D_n$ and let $p = q|_{G''_n}$. Then\\
{\rm a)} the partial $n$-quasigroup $(A, G''_n, p)$ has exactly two extensions;\\
{\rm b)} the $(n-1)$-quasigroups $(A, p^{(0)})$ and $(A, p^{(1)})$ belong to $\overline C_{n-1}$.
\end{lemma}

Therefore,
\[
|\overline C_n \setminus D_n| \leq 2|\overline C_{n-1}|^2.
\]
Moreover, since $p^{(0)}(x) \neq p^{(1)}(x)$ for every $x \in G_{n-1}^1$, we have
\[
|\overline C_n \setminus D_n| \leq \frac{3}{2^n}|\overline C_{n-1}|^2,
\]
and
\[
|\overline C_n| \leq |D_n| + \frac{3}{2^n}|\overline C_{n-1}|^2 \leq 2^{2^n} + \frac{3}{2^n}|\overline C_{n-1}|^2 \quad \text{if } n \geq 6.
\]
Using Proposition~\ref{p3} as the base case, we can easily prove by induction that $|\overline{C}_n| \leq 2 \cdot 2^{2^n}$ provided $n \geq 5$. Since $Q(n) \leq L(n) + |\overline{C}_n|$ by Proposition~\ref{p1}, we obtain
\begin{equation}\label{eq6}
Q(n) \leq (3^{n+1} + 1)2^{2^{n}+1}.
\end{equation}
It follows from \eqref{eq1} and \eqref{eq6} that
\[
Q(n) = 3^{n+1}2^{2^{n}+1} \left(1 + O\!\left(\frac{1}{3^n}\right)\right).
\]

\hrulefill
\renewcommand{\refname}{Further reading (was not a part of the thesis in 2001)}

\end{document}